# ON THE CONCEPT OF (HOMO)MORPHISM : A KEY NOTION IN THE LEARNING OF ABSTRACT ALGEBRA

Thomas Hausberger

Université Montpellier 2

*This article is dedicated to the investigation of difficulties involved in the understanding of the homomorphism concept. It doesn't restrict to group-theory but on the contrary raises the issue of developing teaching strategies aiming at gaining access to structuralist thinking. Emphasis is put on epistemological analysis and its interaction with didactics in an attempt to make Abstract Algebra more accessible.*

## I. INTRODUCTION

In our context, Abstract Algebra means the discipline devoted to the study of algebraic structures, according to the new paradigm established after the publication of van der Waerden's textbook *Moderne Algebra* (1930):

> This image of the discipline turned the conceptual hierarchy of classical algebra upside-down. Groups, fields, rings and other related concepts, appeared now at the main focus of interest, based on the implicit realization that all these concepts are, in fact, instances of a more general, underlying idea: the idea of an algebraic structure. The main task of algebra became, under this view, the elucidation of the properties of each of these structures, and of the relationships among them. Similar questions were now asked about all these concepts, and similar concepts and techniques were used, inasmuch as possible, to deal with those questions. The classical main tasks of algebra became now ancillary. The system of real numbers, the system of rational numbers, and the system of polynomials were studied as particular instances of certain algebraic structures, and what algebra has to say about them depended on what is known about the general structures they are instances of, rather than the other way round (Leo Cory, History of Algebra, Encyclopaedia Britannica Online, 2007).

Abstract algebra is taught in France at third-year University level. The situation reflects the international one and can be summarized by Leron and Dubinsky's (1995, p. 1) statement: *"The teaching of abstract algebra is a disaster, and this remains true almost independently of the quality of the lectures"*. If both students and experienced instructors generally agree upon this provocative claim, this shouldn't prevent us from investigating the nature of the obstacles and proposing alternative or complementary approaches (as Leron and Dubinsky did).

Our approach is based on an epistemological analysis of structuralist thinking (Hausberger 2012) with a view to introducing "meta" aspects in the teaching of abstract algebra. In other words, the goal is to build reflexive activities that can help students to make sense of this particular knowledge. This is justified by the identification, on an epistemological point of view, of the concepts involved as FUGS (formalizing, unifying, generalizing & simplifying) concepts (Robert 1987).

Obstacles regarding the built-up of teaching strategies, using traditional didactic tools such as Brousseau's "théorie des situations" (Brousseau 1986) or Douady's "dialectique outil-objet" (Douady 1986), can be analyzed as consequences of their epistemological nature as FUGS (Rogalski 1995). This approach has been conducted previously by Dorier, Robert, Robinet & Rogalski in the case of linear algebra (Dorier and al. 2000).

In Hausberger 2012, we underlined the fact that several levels of unification have to be distinguished in the context of algebraic structures: at level 1, a general theory applies to quite different objects sharing a common feature (for instance group theory), at level 2, the axiomatic presentation of structures is conducted in a uniform way (leading to structural questions and methods) and puts forward bridges between structures, at level 3, what has been previously a form (structures) is fully taken and studied as object in a superior level of organization (this is category theory or any other mathematical meta-theory of structures, the word "meta" being used here in a different context and with a slightly different meaning. Nevertheless, meta-theories and meta-activities meet in so far as they both introduce a reflexive point of view).

We made the assumption that meta activities of level 2 were worth being introduced to facilitate the access to modern structuralist expositions of algebraic structures. In the case of Linear Algebra, Dorier introduced in fact only meta of level 1. This is justified by the fact that a single structure is at play and that it is the first abstract structure that students encounter and theorize at University.

Moreover, unlike the impression conveyed by Cory in the previous quote, each structure has its own flavor. Indeed, it has its own history of problems, its own "typical objects". There certainly is a tension in modern mathematics related to the articulation of abstract formalism and intuition (requiring a more direct grasp of objects). This tension is acknowledged by philosophers and is also visible in manuals[1]. This raises the following didactic issues and tasks:

−Find a right balance between formalizing and problem-solving involving more concrete objects.

−Think about strategies to gain access to structuralist thinking. This should be progressive: in this respect, the goals set in the teaching of group theory should be different from those regarding ring theory. We will make this statement more precise below.

In our view, the flavor of Linear Algebra resides in interrelating geometry and algebra (enlarging the notion of vector, vector space, or interpreting geometrically reduction theory) and its challenge is also to articulate computational and abstract theoretical aspects (for instance, matrices and linear applications). In group theory, the symmetric group is a paradigmatic example, classifying groups of small order is a very didactical moment and understanding the notion of isomorphism is fundamental to give access to the abstract group concept, thus accomplishing level 1. The emphasis on group action is also fundamental to make the concept fully operational. Ring theory has quite an arithmetical flavor: extending the unique-prime-factorization theorem to rings of algebraic integers has proved to be a motivating force in history for the development of an abstract divisibility theory in which the concept of ideal plays a major role and the main objects to unify, on an elementary level, are numbers and polynomial rings. On the structuralist front, ring theory (being preceded by group theory) should be an opportunity to discuss structural aspects shared in the modern exposition of both theories: the construction of quotients, the concept of homomorphism as structure-preserving function, the isomorphism theorems as tools to compare objects, decomposition theorems into simple objects, characterization by universal properties, etc. These aspects would certainly benefit from being made more explicit through level 2 meta activities. It is nevertheless a real challenge and it remains conjectural whether it is feasible or not at third year University level.

---

1  The preface of Artin 1991 is interesting in this respect and begins with an explicit citation of Hermann Weyl.

This article will now focus on the (homo)morphism and isomorphism concepts which are central in abstract algebra. Students often confess that they lose tracks when homomorphisms come in the foreground after the introduction of isomorphism theorems. Previous didactic studies also report on this issue, for instance Nardi (2000) titles p. 179: "Episode 3: The first isomorphism theorem for groups as a container of compressed conceptual difficulties". We will contribute by deepening the epistemological analysis on morphisms, extending the picture to isomorphism theorems for other structures and conceptions about morphisms in general. This corresponds to our general philosophy of attacking globally the issue of teaching abstract algebra through epistemological, didactic, cognitive studies of processes involved in structuralist thinking.

## II. A FEW CASE STUDIES

### Previous studies on groups

According to Leron, Hazzan & Zazkis (1995), "the very concept of isomorphism is but a formal expression of many general ideas about similarity and differences, most notably, the idea that two things which are different may be viewed as similar under an appropriate act of abstraction". It corresponds to the vague (but crucial and intuitive) idea that two isomorphic groups are essentially the same (on the group-theoretic point of view) and part of the proposed didactic strategy is to help students make sense of it before getting engaged in the formalization of the isomorphism concept. Indeed, such a formalization requires the function concept and reasonable understanding of quantification, which adds further difficulties.

This point of view contrasts with the standard exposition in modern manuals, in which structuralist conceptions have been naturalized. Interestingly, it matches van der Waerden's (1930) who didactically took great care in motivating the introduction of new concepts. If isomorphisms derive formally from homomorphisms, it should be pointed out that isomorphisms come first epistemologically.

Nardi (2000) carries on the work by analyzing tutoring sessions dedicated to retrieving and proving the first isomorphism theorem:

> Let $\phi:G\rightarrow G'$ be a group homomorphism. If $K=ker\ \phi$ then $G/K\simeq Im\phi$. The isomorphism is constructed by setting $\psi(Kg)=\phi(g)$.

She underlines that "the degree of complexity in a problem which requires a well-coordinated [linking $\psi$ and $\phi$] manipulation of mappings between different sets is extremely high" and stresses the numerous "difficulties in the conceptualization of properties associated to the notion of mapping (homomorphic property, 1-1, onto, well-definedness)" as well as the high degree of abstraction involved in the definition of a mapping between the cosets of a subgroup and the elements of the group. Finally, she points out the impact of epistemological arguments that the tutor would put forward to motivate the newly introduced concepts. This certainly brings water to our mill.

**Comparing objects through homomorphisms**

*Context*. The following questions were asked to third year undergraduates as part of a mid-term examination for a second course in abstract algebra devoted to ring and field theory. The student is supposed to retrieve the statement and proof of a classical generalization of the first isomorphism theorem (in the ring context), which has been presented during the lectures as a **tool** for constructing homomorphisms from quotient rings (i.e. of type $A/I\rightarrow B$):

> Let $f:A\rightarrow B$ be a ring homomorphism, $I$ an ideal of $A$ and $\pi:A\rightarrow A/I$. Then $f$ "factorizes through $A/I$" (*i.e.* there exists a homomorphism $f$ such that $f=\bar{f}\circ\pi$) if and only if $I\subset ker\ f$.
>
> The theorem is illustrated by the following commutative diagram:
>
> $$\begin{array}{ccc} A & \xrightarrow{f} & B \\ \pi\downarrow & \nearrow_{\bar{f}} & \\ A/I & & \end{array}$$

Homomorphisms have been introduced as structure-preserving functions, emphasizing the condition $f(1)=1$[2]. To facilitate the retrieval, the ring data has been denoted $(A,+,\cdot,1)$. The general idea has been developed that **homomorphisms aim at "comparing" rings**: if $f:A\rightarrow B$ is bijective then $A$ and $B$ are "essentially the same" and can be identified. If it is only injective (1-1), then $A$ can be identified with a subring of $B$ and the kernel measures the defect of injectivity, if it is surjective (onto), then $B$ can be identified with a quotient of $A$ through the first isomorphism theorem. As an *application* of the generalized theorem, named explicitly "factorization theorem for ring homomorphisms": **Z**$\rightarrow$ **Z**/2**Z** factors

---

[2] which isn't automatic since $A^*$ is not in general a group for the multiplicative law

through **Z**/4**Z** and induces an isomorphism (**Z**/4**Z**)/($\bar{2}$)≃**Z**/2**Z**. The lecture then carried on with the third isomorphism theorem which *generalizes* such isomorphisms and is named "simplification theorem for quotients of quotients". To summarize the didactic intent, the tool-object dialectic (Douady 1986) is at play and epistemological insight was given to connect the formalism with cognitive processes of comparison and identification.

Here are the questions:

1. Recall the factorization theorem for ring homomorphisms and give a proof.

2. Can **Z**/4**Z** be identified with a subgroup of **Z**/8**Z** (in other words, does **Z**/8**Z** contain a subgroup isomorphic to **Z**/4**Z**)? Can **Z**/4**Z** be identified with a subring of **Z**/8**Z** (in other words, does it exist an injective homomorphism **Z**/4**Z**→ **Z**/8**Z**)? Can you construct a ring-homomorphism connecting **Z**/4**Z** and **Z**/8**Z**? Can **Z**/4**Z** be identified with a quotient ring of **Z**/8**Z**?

*Results*. Out of the 13 students who took the test, only 5 managed to retrieve successfully the statement of the theorem (question 1). The diagram was often reproduced, indicating that it played some role in the memorization process. 5 students didn't answer question 1 and among the others, common mistakes concerned the omission of the condition on the ideal, or an inversion of the inclusion. Indeed, the student can have control on this only if he has understood that elements of *I* are mapped to 0 because the diagram commutes. Yet, only 4 students gave an attempt to retrieve the proof: 2 of them successfully demonstrated that the condition is necessary but couldn't carry on so that a single student gave a complete proof. This is a bit surprising since it was emphasized during the lectures that the proof can be worked-out exactly as in the group context. Moreover, the difficulty of articulating $\psi$ and $\phi$ identified by Nardi (2000) is greatly taken in charge didactically by the formulation: "*f* factorizes through the quotient and induces $\bar{f}$". Even after completion of a first course in abstract algebra, the task of elucidating the properties that need to be checked (well-definedness, homomorphic property, 1-1, onto), thus reducing the problem to elementary tasks, is regarded as very complex by students.

A single student gave elements to answer question 2. He began by making the groups more explicit: **Z**/8**Z**={$\bar{0}, \bar{1},..., \bar{7}$}, **Z**/4**Z**={$\bar{0}, \bar{1}, \bar{2}, \bar{3}$} and argues that

G={$\bar{0}, \bar{2}, \bar{4}, \bar{6}$} is a subgroup of **Z**/8**Z** of order 4, therefore isomorphic to **Z**/4**Z**. The map $\phi$: **Z**/4**Z**→ **Z**/8**Z**, $\bar{0}\mapsto\bar{0}, \bar{1}\mapsto\bar{2}, \bar{2}\mapsto\bar{4}, \bar{3}\mapsto\bar{6}$ is an injective morphism of groups but not of rings since $\phi(\bar{1}) \neq \bar{1}$.

Although the student uses the same notation for elements of both quotients, he doesn't mix them and has the right intuition on how to construct the 1-1 group-homomorphism. The teacher would expect a more abstract reasoning: to find a subgroup isomorphic to **Z**/4**Z** amounts to finding an element of order 4 in **Z**/8**Z**. The student doesn't justify that *G* is a group nor that it is cyclic because every subgroup of a cyclic group is also cyclic. It would be interesting to check him on isomorphism classes of groups of order 4. Finally, quantifiers are again handled with too little care: if $\phi$ is not a ring homomorphism, this doesn't imply such a homomorphism doesn't exist. For instance, the student could have argued that **Z**/8**Z** has a unique subgroup of order 4 or used the fact that the $\phi(1)=1$ property completely determines $\phi$ (which then isn't a homomorphism). The last part of the question remains unanswered and indicates that finding a map between two quotient rings is a conceptually hard task (and a bit tricky one since the arrow is reversed: **Z**/4**Z**←**Z**/8**Z**), even if it consists in reproducing the reasoning made during the course on a very close example.

In fact, the two questions were asked before heavy use of the factorization theorem for rings was made during the tutoring sessions. The experiment shows that the test was premature even if the course was epistemologically-oriented and the students could rely on previously acquired conceptual knowledge (which remains very fragile). We will engage in further investigations of the difficulties in 2013 by asking the same questions during tutoring work. Nevertheless, this experiment already confirms that the homomorphism concept is much of an obstacle, even after the completion of a first course in abstract algebra devoted to group theory.

## III. MORE EPISTEMOLOGICAL INSIGHT

**By the way, what is a structure ?**

Abstract Algebra teachers often speak about structures... but they never mathematically define any concept of structure! Let us explain this a priori abnormal phenomenon by giving a few details concerning the didactic transposition (Chevallard 1985) of the notion of structure.

In his attempt to give an historical account, Cory (1996) makes the distinction between "body of knowledge" and "image of knowledge". Interestingly, the notion of structure takes its origin in the latter:

> This textbook [*Moderne Algebra*] put forward a new image of the discipline that implied in itself a striking innovation: the structural image of algebra. In the forthcoming account, it is this specific, historically conditioned image of mathematical knowledge that will be considered as implicitly defining the idea of a mathematical structure (Cory 1996, p.8).

We won't comment on historical methodology, but it is certainly the combination of Noether's mathematical ideas (see below) and a didactic intent to expose the recent advances in algebra in a systematic and clear fashion, in an organized and integrated whole, that lead to the idea of structure. It remained implicit in so far as van der Waerden didn't give any comment, formal or non-formal, on what he meant by a 'structure'. Bourbaki, on the contrary, gave a formal-axiomatic elucidation of the concept of mathematical structure in the first book of his treatise '*Elements de mathématiques*' dedicated to set theory and published in the 1950s. In parallel, he promoted (Bourbaki 1948) the structural image of mathematics within the noosphere (Chevallard 1985). Yet, Bourbaki's definition did hardly play any role in the exposition: it only provided a general framework which in fact didn't prove to be mathematically functional (Cory 1996, p. 324, see also Mac Lane 1996), unlike category theory which is very advanced and a too hovering viewpoint for the present purpose. Therefore no definition at all is given in more recent manuals.

As a consequence, students are supposed to learn by themselves and by the examples what is meant by a structure whereas sentences like "a homomorphism is a structure-preserving function" is supposed to help them make sense of a homormorphism. Is that possible without any clarification on the notion of structure? The students certainly understand that when we talk about structures we refer to sets of axioms and we say that a mathematical object has a particular structure when these are fulfilled. This reflects the fact that the notion of structure is an outgrowth of the widespread use of the axiomatic method[3]. In our view, meta-discourse on axiomatics together with activities devoted to building axiomatics and "playing" with the axioms would help students to make sense of them and memorize them. It is not straightforward to comprehend that axiomatics encode properties of mathematical objects which are analyzed abstractly as being made of elements connected by relations (therefore primitive terms of axiomatics are sets equipped with extra data encoding the relations such as laws of composition or binary relations). The abstract formalism (axioms and the language of set theory) is much of an obstacle and hides the simple ideas. Too little emphasis is made on the idea of relations: structure-preserving is synonymous with operation-preserving in contemporary manuals, unlike van der Waerden's (1930), in which examples are given of non-algebraic structures such as ordered sets when introducing the homomorphic property, defined as preserving relations. The idea of relation is also useful to understand quotients: making a quotient is equivalent to introducing more relations. This should help the students to interpret, for instance, the quotient $\mathbf{Z}[X]/(10X-1)$ as a ring isomophic to the decimal numbers. We will need of course to engage in further didactical studies in order to support our conjectural claims regarding conceptions and teaching strategies.

Coming back to the epistemological investigation, the word 'structure' is used in fact in two more contexts with a different meaning, which may induce some confusion: we also want to identify the different isomorphism classes for a given structure and we say for instance that the abstract group-structure of $\mathbf{Z}/3\mathbf{Z}$ is that of a cyclic group of order 3. Finally, we call 'structure-theorem' a result describing the way an object can be reconstructed from simpler objects of the same type.

---

3  Historically, this method was initially deployed by Hilbert to remedy the imperfections of Euclid's axioms for geometry and give rigorous descriptions of mathematical objects suitable for applying the demonstrative apparatus. In this trend, logical use of the method is about questioning consistency, mutual independence, completeness of the set of axioms, in order to give foundations to the theory. In our context, we are concerned with a more immanent use of axiomatics: properties of concrete objects are abstracted in order to get a better understanding of ideas underlying mathematical constructions and proofs. This allows a unification of objects and methods, leading to the structural point of view.

**Noether's set-theoretic foundation of Algebra**

We have already distinguished the concept of isomorphism as an equivalence relation from that of a function (with properties) which relies on the notion of homomorphism. The equivalence relation is in fact a crucial concept with regard to the process of abstraction, which often amounts to selecting common characteristics of objects being thus taken as equivalent. A classification is a description of the corresponding partition: for instance, group theory will classify groups of a given order. We will now focus on the latter conception of isomorphism: it is the heritage of Noether who developed

> […] what she called her set-theoretic foundations for algebra. This was not what we now call set theory. […] Rather, her project was to get abstract algebra away from thinking about operations on elements, such as addition or multiplication of elements in groups or rings. Her algebra would describe structures **in terms of selected subsets** (such as normal subgroups of groups) **and homomorphisms** (MacLarty, 2006, p. 188).

We pretend that Noether's new conceptual approach (characterized in the last sentence of the quote) is a **major epistemological difficulty** in the learning of abstract algebra. It is also the key to level-2 structuralism (referring to the introduction): this indeed allows a unified treatment of structures. It proved in history to be a major breakthrough leading to a complete rewriting of algebra (in terms of newly-forged concepts that emerged from the new methodology: noetherian rings, principle-ideal domains,...).

Dedekind's theory of ideals in which divisibility relation between (algebraic) integers were replaced by inclusion of ideals certainly contributed to the transition from an arithmetical conception of algebra to the set-theoretic conception. But Noether's chief tools were isomorphism theorems and she made it obvious that this applied for different kinds of structures. She considered only onto homomorphisms, denoted $M \sim \bar{M}$ (the functional notation f:$M \rightarrow \bar{M}$ comes from topology) and correlated them to distinguished classes of subsets through the first isomorphism theorem. As an illustration of the generality of the principle:

> Ideals bear the same relation to ring homomorphisms as do normal subgroups to group homomorphisms. Let us start from the notion of homomorphism (van der Waerden, 1949, p. 51).

**Homomorphisms or morphisms?**

On forums[4] are taking place interesting discussions on the differences to be made between homomorphisms and morphisms. Some people argue that homomorphism is the old terminology and the shorter word should be adopted for pragmatic reasons. Others mention that morphisms come from category theory (which remains obscure to them). Finally, one of them argues that the difference between homomorphism and homeomorphism is clear, but not between homomorphism and morphism. This leads to a confrontation of two different definitions of a morphism in the context of topological spaces: open maps (preserving open sets) versus continuous maps. But this didn't allow to make the point underneath the morphism concept which derives from category theory and thus the difference to be made.

---

4  For instance, on the site http://www.les-mathematiques.net. We unfortunately have no space for a transcript here.

Indeed, the morphism concept is a relativization of the homomorphism concept as a structure-preserving function (which doesn't apply to morphisms in topology since they preserve open sets by *inverse* image): in a category, one is free to decide which maps are morphisms, these define the category together with a given type of objects. If morphisms tend to replace homomorphisms, beyond the pragmatic argument, it might be that algebra practitioners are implicitly assuming that they are working in a category.

## IV. CONCLUDING REMARKS AND PERSPECTIVES

This study certainly contributes to break the illusion of transparency concerning the concept of homomorphism and the idea of a mathematical structure. It aims at reestablishing the rationale of this particular knowledge through the epistemological investigation of the concept and the engineering of epistemologically-oriented activities. In this spirit, a simple didactical situation has been given in part II. Unfortunately, it needs to be re-experimented in order to reveal its full potential. Meanwhile, we have engineered, on the basis of the epistemological analysis of the notion of structure presented in this paper, an activity dedicated to reflecting on the axiomatic method and the structuralist viewpoint in a simple context ("mini-theory"). We hope to report soon on the results of the pre-experimentation of this activity. On a more theoretical perspective, following Chevallard's idea that a "body of knowledge" is a complex of praxeologies (see for instance Chevallard 2002), we will need to track and analyze the "structuralist praxeologies" in order to support the idea, presented in part I of this paper, that there really is a "structuralist agenda" to meet in relation with the learning of abstract algebra. Finally, our analysis connects the transition problem that occurs at third year university level in relation to abstract algebra with epistemological transitions: the systematization of the axiomatic method, after Hilbert, and the transition, after Noether, from thinking about operations on elements to thinking in terms of selected subsets and homomorphisms. This analysis should be crossed with a didactical investigation of the transition problem: the analysis of praxeologies should again bring much light.